\documentclass[12pt,a4paper,reqno]{amsart}

\renewcommand{\bar}{\overline}

\newcommand{\pa}{\partial}
\renewcommand{\phi}{\varphi}
\newcommand{\mi}{{\mathcal M}}

\newcounter{hours}\newcounter{minutes}

\newcommand{\ka}{K\"ahler }

\font\strange=msbm10

\newcommand{\R}{{{\mathchoice  {\hbox{$\textstyle{\text{\strange R}}$}}
{\hbox{$\textstyle{\text{\strange R}}$}}
{\hbox{$\scriptstyle  N\kern-0.3em  R$}}  
{\hbox{$\scriptscriptstyle  R\kern-0.2em  R$}}}}}

\newcommand{\Z}{{{\mathchoice  {\hbox{$\textstyle{\text{\strange Z}}$}}
{\hbox{$\textstyle{\text{\strange Z}}$}}
{\hbox{$\scriptstyle  Z\kern-0.3em  Z$}}
{\hbox{$\scriptscriptstyle  Z\kern-0.2em  Z$}}}}}

\newcommand{\N}{{{\mathchoice  {\hbox{$\textstyle{\text{\strange N}}$}}
{\hbox{$\textstyle{\text{\strange N}}$}}
{\hbox{$\scriptstyle  N\kern-0.3em  N$}}
{\hbox{$\scriptscriptstyle  N\kern-0.2em  N$}}}}}

\newcommand{\frk}[1]{{\mathfrak{#1}}}

\newcommand{\bb}{{\frac{\sqrt{-1}}{2\pi}}}

\usepackage{amsmath,amsthm,amscd}

\usepackage{calc}

\title 
{Gradient estimates of the Yukawa coupling}
\author{Zhiqin Lu}
\date{September 11, 2001}
\subjclass{Primary: 53A30; Secondary: 32C16}

\keywords{Gradient estimates, Calabi-Yau threefold, 
Yukawa coupling}
\address[Zhiqin Lu]
{Department of Mathematics\\
University of California, Irvine\\ 
Irvine, CA 92697}
\email[Zhiqin Lu]{zlu@math.uci.edu}
\thanks{Research supported by NSF grant DMS 0196086}

\newtheorem{theorem}{Theorem}

\newtheorem{lemma}{Lemma}

\newtheorem{prop}{Proposition}

\newtheorem{definition}{Definition}

\newtheorem{example}{Example}
\theoremstyle{remark}
\newtheorem{rem}{Remark}
\newtheorem*{co}{\bf Conjecture}

\begin{document}
\leftline{{\footnotesize   Special Geometric Structures in String Theory}}
\leftline{{\footnotesize  Bonn, 8th--11th September, 2001}}
\vskip 1.5 true cm

\maketitle


\tableofcontents

\section{Introduction}

A polarized Calabi-Yau manifold is a pair $(X,\omega)$ of a compact
algebraic manifold  $X$ with zero first Chern class and a K\"ahler form
$\omega\in H^2(X,\mathbb Z)$. The form $\omega$ is called
a  polarization. Let
${\mathcal M}$ be the universal deformation space of $(X,\omega)$.
${\mathcal M}$ is smooth by the  theorem of Tian~\cite{T1}.
By~\cite{Y1}, we may assume that each $X'\in{\mathcal M}$ is a
K\"ahler-Einstein manifold. i.e.
the associated K\"ahler metric $(g'_{\alpha\bar\beta})$ is Ricci flat. The
tangent space $T_{X'}{\mathcal M}$ of ${\mathcal M}$ at $X'$ can be
identified with $H^1(X',T_{X'})_{\omega}$ where
\[
H^1(X',T_{X'})_{\omega}=\{\phi\in H^1(X',T_{X'})|
 \phi\lrcorner\omega=0\}.
\]

The Weil-Petersson metric
$G_{WP}$ on ${\mathcal M}$ is defined by
\[
G_{WP}(\phi,\psi)=\int_{X'} 
{g'}^{\alpha\bar\beta}g'_{\gamma\bar\delta}\phi^\gamma_{\bar\beta}\bar{\psi}^\delta_{\bar\alpha}
dV_{g'},
\]
where $\phi=\phi^\gamma_{\bar\beta}\frac{\partial}{\partial z^\gamma}
d\bar{z}^\beta$, 
$\psi=\psi^\delta_{\bar\alpha}\frac{\partial}{\partial z^\delta}
d\bar{z}^\alpha$ are in $H^1(X',T_{X'})_\omega$, $g'=g'_{\alpha\bar\beta}
dz^\alpha d\bar{z}^\beta$ is the K\"ahler-Einstein metric on $X'$
associated with the polarization $\omega$.

A natural question on $\mi$ is that whether the Weil-Petersson
metric is complete. In~\cite{Lu8}, the author proved that there are
no non-trivial complete special \ka manifolds. 
If $\dim M=3$, then
$\mi$ is a
projective special \ka manifold. The corresponding
conjecture in this case would be:  

\begin{co}  
If the moduli space $\mi$ of a Calabi-Yau threefold is complete with
respect to the Weil-Petersson metric, then it is locally
symmetric.
\end{co}

\begin{rem} 
The author learned from the referee that homogeneous 
projective special \ka manifolds of semisimple group 
were classified by Alekseevsky and Cort\'es~\cite{AC}.
They are all Hermitian symmetric of noncompact type.
\end{rem}

The list of the homogeneous projective manifolds are:

\begin{enumerate}
\item[A).] ${\frk l}=\frk{sl}_{n+2}(\mathbb C), \quad
M=\mathbb C\mathbb H^{n-1}=SU_{1,n-1}/S(U_1\cdot
U_{n-1})$;
\item[BD).] $\frk l=\frk{so}_{n+4}{\mathbb C},
\quad M=(SL_2(\mathbb R)/SO_2)\times
(SO_{2,n-2}/SO_2\cdot SO_{n-2})$;
\item[G).] $\frk l=\frk g_2(\mathbb C), \quad
M=\mathbb{CH}^1=SL_2(\mathbb R)/SO_2$;
\item[F).] $\frk l=\frk f_4(\mathbb C), \quad 
M=Sp_3(\mathbb R)/U_3$;
\item[E6).] $\frk l=\frk e_6(\mathbb C), \quad
M=SU_{3,3}/S(U_3\cdot U_3)$;
\item[E7).] $\frk l=\frk e_7(\mathbb C), \quad
M=SO^*_{12}/U_6$;
\item[E8).] $\frk l=\frk e_8(\mathbb C), \quad
M=E_7^{(-25)}/E_6\cdot SO_2$.
\end{enumerate}

In \S 3, we give an example of the complete locally
symmetric space, which coresponds to the above type (A),
the complex hyperbolic space.

For the rest of this paper, we will
concentrate on  moduli space of Calabi-Yau threefolds.
General moduli space will be considered
elsewhere~\cite{Lu20}. 

In the three dimensional case,
associated to the Weil-Petersson metric is the Yukawa coupling. It
can be defined as 
\[
F(\phi,\psi,\xi)=\int_{X'}\phi\wedge
\psi\wedge\xi\lrcorner\Omega,
\]
where $\Omega$ is a $(3,0)$ form on $X'$. Since $X'$ is Calabi-Yau,
$(3,0)$ forms on $X'$ differ by  constants.

Note that the Yukawa coupling depends not only on 
$\phi, \psi,\xi$ but depends on $\Omega$ as well. In fact, one can
prove ~\cite{Gr} that
$F$  is a holomorphic section of the bundle
${\rm Sym}^3(T^*M)\otimes (\underline F^3)^{\otimes 2}$, where
$\underline F^3$ is the first Hodge bundle on $\mi$(cf. ~\cite{Gr}).

One of the fundamental properties of the Weil-Petersson metric and
the Yukawa coupling is that  they can be defined 
``extrinsically'' in the sense that they can be defined 
only using the fact that the moduli space is a horizontal slice. 
In fact, let $Q(\Omega,\bar\Omega)$ be defined in ~\eqref{qq}.
We
have
\[
\omega_{WP}=-\pa\bar\pa\log Q(\Omega,\bar\Omega),
\]
(cf. ~\cite{T1}), 
where $\omega_{WP}=\bb h_{i\bar j} dz_i\wedge d\bar z_j$ is the \ka
form of the Weil-Petersson metric, and
\[
F_{ijk}=Q(\pa_i\pa_j\pa_k\Omega, \Omega).
\]

One can contract the Yukawa coupling to get the following 
$(1,1)$ tensor  
\begin{equation}\label{PP}
P=\bb P_{i\bar j}dz^i\wedge d\bar z^j= \bb h^{p\bar q} h^{r\bar s}
  F_{ipr}\bar{F_{jqs}} dz^i\wedge d\bar z^j.
\end{equation}
This tensor is important because 
of the following theorem~\cite{Lu5}:

\begin{theorem}\label{thm11}
Let $\omega_H=2\omega_{WP}+P$, and let $n=\dim\,\mi$. Then
\begin{enumerate}
\item $\omega_H$ is a K\"ahler metric on ${\mathcal M}$;
\item The holomorphic bisectional curvature of $\omega_H$ is 
nonpositive.
Furthermore,
Let $\alpha=((\sqrt{n}+1)^2+1)^{-1}>0$. Then the  Ricci curvature 
${\rm Ric}(\omega_H)\leq
-\alpha \omega_H$ and the holomorphic sectional curvature 
is also less
than or
equal to $-\alpha$.
\item If ${\rm Ric}(\omega_H)$ is bounded, then the Riemannian
sectional curvature of $\omega_H$ is also bounded.
\end{enumerate}
\end{theorem}
We call $\omega_H$ the Hodge metric on $\mi$.

If the Yukawa coupling is bounded, then the Weil-Petersson
metric and the Hodge metric are equivalent. The other side of the
theorem is the main result of this paper:

\begin{theorem}\label{thm3}
Assume that the Weil-Petersson metric is complete. Then
there is a constant $C_1(m,n)$, depending only on $m,n$, such that
\[
|\nabla^m F|^2\leq C_1(m,n),
\]
for any nonnegative integer $m$, where $\nabla$ is the Hermitian 
connection of the bundle ${\rm Sym}^3(T^*M)\otimes (\underline
F^3)^{\otimes 2}$ and $n$ is the complex dimension of the moduli
space
$\mi$.
\end{theorem}

\noindent
{\bf Acknowledgment.} This paper is a refinement of  the talk the
author gave at the workshop ``Special Geometric Structures and
String Theory'' at Universit\"at  Bonn in September 2001. The author
thanks the organizers Dmitri V. Alekseevsky, Vicente Cort\'es,
Chand Devchand and Antoine Van Proeyen for the invitation and helpful
conversations during the workshop.

\section{On the Complete Weil-Petersson Metric}
In this section, we use the method of gradient estimate to prove
Theorem~\ref{thm3}, the main result of this paper. We first prove
the following weak version of Theorem~\ref{thm3}.

\begin{theorem}\label{thm21}
Suppose ${\mathcal M}$ is the moduli space of a Calabi-Yau
threefold. If the Weil Petersson metric 
on ${\mathcal M}$ is complete, then the norm of the
Yukawa coupling   with respect to the Weil-Petersson
metric is bounded. On the other hand, if the Yukawa
coupling with respect to the Weil-Petersson metric is
bounded, then the Weil-Petersson metric is equivalent to
the Hodge metric.
\end{theorem}
  
\noindent
{\bf Proof:} Let $\omega_H$ be the Hodge  metric, 
then by definition
\begin{equation}
\omega_H=2\omega_{WP}+P,
\end{equation}
where $P$ is defined in~\eqref{PP} and
$\omega_{WP}=\bb h_{i\bar j} dz_i\wedge d\bar z_j$ is the \ka form
of the Weil-Petersson metric. The trace of the tensor $P$ is the
norm of the Yukawa coupling. Thus the Hodge
metric and the Weil-Petersson metric are mutually equivalent
if the Yukawa coupling is bounded.

On the other side,
let
\[
f=|F|^2=e^{2K}F_{ijk}\bar{F_{abc}}h^{i\bar{a}}
h^{j\bar{b}}h^{k\bar{c}},
\]
where $K=-log\sqrt{-1}Q(\Omega,\bar\Omega)$, and 
\begin{equation}\label{qq}
Q(\Omega,\bar\Omega)=\int_{X'}\Omega\wedge\bar\Omega.
\end{equation}
We have
\begin{align*}
&f_{\bar\alpha}=2e^{2K}K_{\bar\alpha}F_{ijk}\bar{F_{abc}}
h^{i\bar a}h^{j\bar b}h^{k\bar c}\\
&+e^{2K}F_{ijk}\bar{\pa_\alpha F_{abc}}
h^{i\bar a}h^{j\bar b}h^{k\bar c}
+e^{2K}F_{ijk}\bar{F_{abc}}\pa_{\bar\alpha}(
h^{i\bar a}h^{j\bar b}h^{k\bar c})
\end{align*}
for $\alpha=1,\cdots,n$.
Thus we have
\[
\Delta f=e^{2K}|\pa_\alpha F_{ijk}+2K_\alpha F_{ijk}|^2
+2n|F|^2
+e^{2K}F_{ijk}\bar{F_{abc}}\pa_\alpha\bar{\pa_\alpha}
(h^{i\bar a}h^{j\bar b}h^{k\bar c}),
\]
where $\Delta$ is the complex Laplacian of the Weil-Petersson
metric. Under the normal coordinates,
\[
e^{2K}F_{ijk}\bar{F_{abc}}\pa_\alpha\bar{\pa_\alpha}
(h^{i\bar a}h^{j\bar b}h^{k\bar c})
=3e^{2K}F_{ijk}\bar{F_{ajk}}R_{a\bar i}\]
where $R_{a\bar i}$ is the Ricci curvature of the Weil-Petersson
metric. Thus
\[
\Delta f\geq 2n|F|^2+3e^{2K}F_{ijk}\bar{F_{ajk}}R_{a\bar i}.
\]
It is known in~\cite{S} that
\begin{equation}
R_{a\bar i}=-(n+1)\delta_{a\bar i}+e^{2K}F_{amn}\bar{F_{imn}}.
\label{eq:eq3}
\end{equation}
Thus
\begin{align*}
&\Delta f \geq 2n|F|^2-3(n+1)|F|^2
+3e^{4K}\sum_{a,i}|\sum_{j,k}F_{ijk}\bar{F_{ajk}}|^2\\
&\geq -(n+3)|F|^2+3e^{4K}\sum_i(\sum_{j,k}|F_{ijk}|^2)^2\\
&\geq -(n+3)|F|^2+\frac 3n e^{4K}(\sum_{i,j,k}|F_{ijk}|^2)^2\\
&=\frac 3n f^2-(n+3)f.
\end{align*}

We now recall a version of the maximum principle from~\cite{T1}.

\begin{prop}\label{prop1}
Suppose that $(M,g)$ is a complete K\"ahler manifold. If the Ricci
curvature of $g$ is bounded from below and $\varphi$ is a nonnegative
function satisfying
\[
\Delta\varphi\geq c_1\varphi^\alpha-c_2\varphi-c_3,
\]
where $\alpha>1, c_1>0, c_2,c_3\geq 0$ are constants. then
\[
\sup\varphi\leq Max (1,(\frac{c_2+c_3}{c_1})^{\frac{1}{\alpha}}).
\]
\end{prop}

By equation~\eqref{eq:eq3} we know that the Ricci curvature is
bounded from
below. Thus using Proposition~\ref{prop1}, we have
\[
f\leq\sqrt{\frac{n(n+3)}{3}}.
\]
\qed
\begin{rem}
We can also get similar estimates on moduli spaces with incomplete
Weil-Petersson metric. In that case, a different version of Maximum
principle should be set up.
\end{rem}

\noindent
{\bf Proof of Theorem~\ref{thm3}.} We define
\begin{equation}\label{6-1}
 f_m=|\nabla^m F|^2
\end{equation}
for $m=0,1,2,\cdots$.
The inequality is true for $m=0$ by Theorem~\ref{thm21}. Assume
that the inequality is also true for all $0\leq i\leq m-1$. That is,
we have a constant $\tilde C_1(m,n)$ such that
\begin{equation}\label{6}
|\nabla ^i F|^2\leq\tilde C_1(m,n)
\end{equation}
for any $0\leq i\leq m-1$. We are going
to prove that $|\nabla ^m F|$ is bounded. First we have the
following lemma:

\begin{lemma}\label{lem1}
With the above assumption, there is a constant $C_2$ depending
only on $m,n$ such that
\[
\Delta f_m\geq f_{m+1}-C_2(m,n) (f_m+1).
\]
\end{lemma}
\noindent
{\bf Proof.} By~\eqref{6-1}, we have
\begin{align}\label{7}
\begin{split}
&\Delta f_m=|\nabla^{m+1} F|^2+|\bar\pa\nabla^m F|^2\\
&+<h^{i\bar j} \nabla_i\bar\nabla_j\nabla^m F,\bar
{\nabla^m F}>+<\nabla^m F,\bar{h^{j\bar
i}\bar\nabla_i\nabla_j\nabla^mF}>.
\end{split}
\end{align}
Changing the order of the covariant derivative, we get
\begin{equation}\label{8}
\bar\nabla_i\nabla_j\nabla^m F=R(\frac{\pa}{\pa\bar z^i},
\frac{\pa}{\pa z^j}) \nabla^m F+\nabla_j\bar\nabla_i\nabla^m F.
\end{equation}
By the Strominger's formula~\cite{S} of the curvature of the
Weil-Petersson metric
\begin{equation}\label{9}
R_{i\bar jk\bar l}=h_{i\bar j}h_{k\bar l}
+h_{i\bar l}h_{k\bar j}-e^{2K} h^{p\bar q} F_{ikp}\bar{F_{jlq}},
\end{equation}
and by the assumption~\eqref{6}, we see that
\begin{equation}\label{10}
|\nabla^m R_{i\bar jk\bar l}|\leq C_2 |\nabla^m F|+C_3
\end{equation}
for some constants $C_2$ and $C_3$ depending only on $m,n$. Thus
by~\eqref{8}
\begin{equation}\label{10-1}
|h^{i\bar j} \bar\nabla_i\nabla_j\nabla^m F|+|h^{i\bar j}
\nabla_i\bar\nabla_j\nabla^m F|\leq 
C_4 f_m+C_5
\end{equation}
for constants $C_4$ and $C_5$ depending only on $m,n$. By~\eqref{7}
 and~\eqref{10-1}, there is a constant $C_2(m,n)$ such 
that
\[
\Delta f_m\geq f_{m+1}-C_2(m,n) (f_m+1).
\]
In particular, using ~\eqref{6}, we have
\begin{equation}\label{10-2}
\Delta f_i\geq f_{i+1}-C_4
\end{equation}
for $0\leq i\leq m-1$, and the constant $C_4$ depending only on
$m,n$.

\qed

\noindent
{\bf Continuation of the proof of Theorem~\ref{thm3}.}
It is not hard to see that
\begin{equation}\label{11}
|\nabla f_m|\leq 2\sqrt{f_{m+1} f_m}.
\end{equation}
Let 
\[
g_m=f_m(A+f_{m-1}),
\]
where constant $A$ is to be determined. Then using Lemma~\ref{lem1}
and~\eqref{10-2},\newline~\eqref{11}, we have
\begin{equation}\label{12}
\Delta g_m\geq A f_{m+1}+f_m^2-C_5(A) f_m-C_6(A)-C_7
f_m\sqrt{f_{m+1}},
\end{equation}
where $C_5(A)$, $C_6(A)$ are constants 
depending only on $m,n$ and $A$ and $C_7$ is the constant depending
only on $m,n$. We choose that $A=C_7^2$. Then
\[
Af_{m+1}+\frac 14 f_m^2\geq C_7 f_m\sqrt{f_{m+1}}.
\]
Thus
there are constants $\delta>0$ and $C_{8}$, depending
only on $m,n$, such that
\begin{equation}\label{13}
\Delta g_m\geq \delta g_m^2-C_{8}.
\end{equation}
Using  the maximal principal Proposition~\ref{prop1}, 
\[
g_m\leq C_8/\delta+1.
\]
Since we may have chosen $A\geq 1$, we have
\[
f_m\leq g_m\leq C_8/\delta+1,
\]
and the theorem is proved.
\qed

\section{An example}
In this section, we give an example of locally symmetric
horizontal slice. In~\cite{AC}, a complete list of 
homogeneous projective special \ka manifolds is given.

We first introduce the notion of
classifying space by recalling the definitions and notations
in~\cite{Gr}.

Suppose $X$ is a simply connected algebraic Calabi-Yau three-fold. The
Hodge
decomposition of the cohomology group $H=H^3(X,C)$ is
\[
H^3(X,C)=H^{3,0}\oplus H^{2,1}\oplus H^{1,2}\oplus H^{0,3},
\]
where
\[
H^{p,q}=H^q(X,\Omega^p),
\]
and $\Omega^p$ is the sheaf of the holomorphic $p$-forms.
The  skew-symmetric form $Q$ on $H$  is defined by
\[
Q(\xi,\eta)=-\int_X\xi\wedge\eta.
\]
By the Serre duality and the fact that the canonical bundle is
trivial,
$\dim H^{2,1}=\dim H^{1,2}=\dim H^1(X,T_X)=n$, and 
$\dim H^{3,0}
=\dim H^{0,3}=1$. Thus $H^3(X,C)=C^{2n+2}$ is a (2n+2)-dimensional
complex vector
space.

It is easy to check that $Q$ is skew-symmetric. Furthermore, we have the
following two Hodge-Riemannian relations:

1. $Q(H^{p,q}, H^{p',q'})=0$ unless $p'=3-p$ and $q'=3-q$;

2. $(\sqrt{-1})^{p-q}Q(\psi,\bar\psi)>0$ for any nonzero element
$\psi\in H^{p,q}$.

We define the Weil operator $C:H\rightarrow H$ by
\[
C|_{H^{p,q}}=(\sqrt{-1})^{p-q}.
\]
For any collection of $\{H^{p,q}\}$'s, set
\begin{align*}
& F^3=H^{3,0};\\
& F^2=H^{3,0}\oplus H^{2,1};\\
& F^1=H^{3,0}\oplus H^{2,1}\oplus H^{1,2}.
\end{align*}
Then $F^1,F^2,F^3$ defines a filtration of $H$
\[
0\subset F^3\subset F^2\subset F^1\subset H.
\]
Under this terminology, the Hodge-Riemannian relations can be
re-written as

3. $Q(F^3,F^1)=0, Q(F^2,F^2)=0$;

4. $Q(C\psi,\bar\psi)>0$ if $\psi\neq 0$.

\vspace{0.2in}

Now we suppose that $\{h^{p,q}\}$ is a collection of integers such that
$p+q=3$ and $\sum h^{p,q}=2n+2$. 

\begin{definition}
With the notations as above, the classifying space $D$ of the Calabi-Yau
three-fold is the set of all collection of subspaces $\{ H^{p,q}\}$ of $H$
such that 
\[
H=\underset{p+q=3}{\oplus}H^{p,q}\qquad
H^{p,q}=\bar{H^{q,p}},\qquad
 \dim\,H^{p,q}=h^{p,q},
\]
and on which $Q$ satisfies the two Hodge-Riemannian relations 1,2.

Set $f^p=h^{n,0}+\cdots +h^{p,n-p}$. Then $D$ is also
the
set of all filtrations
\[
0\subset F^3\subset F^2\subset F^1\subset H, \qquad 
F^p\oplus\bar{F^{4-p}}=H
\]
with $\dim F^p=f^p$ 
on which $Q$ satisfies the bilinear relations 3,4.
\end{definition}

$D$ is a homogeneous complex manifold. 
The horizontal distribution $T_h(D)$ is defined as
\[
T_h(D)=\{X\in T(D)| XF^3\subset F^2,XF^2\subset F^1\},
\]
where $T(D)$ is the holomorphic tangent bundle which can be identified as a
subbundle of the 
(locally trivial) bundle 
$Hom(H^3(X,C),H^3(X,C))$. So $X$ naturally acts on $F^p$.

\begin{definition}
A complex integral submanifold of the horizontal distribution $T_h(D)$  is
called a horizontal slice. 
\end{definition}

Suppose $U\subset {\mathcal M}$ is a neighborhood of ${\mathcal M}$ at the
point $X$. Then there is a natural map $p: U\rightarrow D$, called the
period map, which sends a Calabi-Yau threefold to its ``Hodge Structure''.
To be precise, Let $X'\in U$. Then there is a natural identification of
$H^3(X',C)$ to $H^3(X,C)=H$. So $\{H^{p,q}(X')\}_{p+q=3}$ are the
subspaces of $H$ satisfying the Hodge-Riemannian Relations. We define
$p(X')=\{H^{p,q}(X')\}\in D$. It is proved in~\cite{Gr} that
$p(U)$ is a horizontal slice.

Now we introduce a result of  Bryant and Griffiths~\cite{BG}. Their
results can
be
briefly
written as follows: 

We assume that $eV\in U$. i.e. the  horizontal  slice
passes the original point of $D$, where the original point is defined
as $\{f^3,f^2,f^1\}\in D$ as follows:
there is a basis $e_1,\cdots,e_{2n+2}$ of $H$ under which $Q$
can be represented as
\[
Q=
\sqrt{-1}
\left(
\begin{array}{cc}
& 1\\
-1 &
\end{array}
\right).
\]
If we let
\begin{align*}
&f^3=span \{e_1-\sqrt{-1}e_{n+2}\},\\
&f^2=span\{e_1-\sqrt{-1}e_{n+2},e_2+\sqrt{-1}e_{n+3},
\cdots, e_{n+1}+\sqrt{-1}e_{2n+2}\},
\end{align*}
and $f^1$ is the hyperplane perpendicular to $f^3$ with respect
to $Q$, then
\[
\{0\subset f^3\subset f^2\subset f^1\subset H\}\in D.
\]

According to Bryant and Griffiths, there is a holomorphic
function $u$ with $u(0)=-\sqrt{-1}$, 
$\nabla u(0)=0$ and $\nabla^2u(0)=\sqrt{-1} I$ ($I$ is the
identity matrix) defined on
a neighborhood of the original point of 
${\mathbb C}^n$ such
if $(z^1,\cdots, z^n)$ is the local holomorphic coordinate of $U$ at $eV$,
the original point, then the horizontal slice passing through
$eV$ can be represented by
\begin{equation}\label{o}
F^3=span(1,\frac{1}{\sqrt{2}}z_1,\cdots,\frac{1}{\sqrt{2}}z_n,
u-\sum_i\frac
12z_iu_i,\frac{1}{\sqrt{2}}u_1,\cdots,\frac{1}{\sqrt{2}}u_n),
\end{equation}
and   $F^2=\nabla F^3$,  
$F^1\perp F^3$ via $Q$. 

\begin{example} With the above notations, 
we choose
\[
u=-\sqrt{-1}+\frac{\sqrt{-1}}{2} \sum_{i=1}^nz_i^2.
\]
Then we can define the horizontal slice with $F^3$ being given by
\[
F^3=(1,\frac{1}{\sqrt{2}}z_1,\cdots,\frac{1}{\sqrt{2}}z_n,
-\sqrt{-1}, 
\frac{\sqrt{-1}}{\sqrt 2}z_1,
\cdots, \frac{\sqrt{-1}}{\sqrt 2}z_n).
\]
\end{example}

The Yukawa coupling of the above example is identically zero. Thus
by~\eqref{9}, the curvature tensor is
\[
R_{i\bar jk\bar l}=h_{i\bar j}h_{k\bar l}+h_{i\bar l}h_{k\bar j},
\]
which  is parallel. In order to see that the horizontal slice
we defined is complete, we first observed that since the Yukawa
coupling is zero, the Hodge metric is two times the Weil-Petersson
metric because $P\equiv 0$ in~\eqref{PP}. Using~\cite[Lemma
3.8]{Lu5}, we can isometrically embed the horizontal slice to the
Siegel manifold
$H_{n+1}$, where
$H_{n+1}$ is the set of all $(n+1)\times(n+1)$ matrices of the form
$X+\sqrt{-1}Y$ with $X, Y$ symmetric and $Y$ positive.
By~\cite{Lu5}, the embedding can be represented by
the matrix
\[
\begin{pmatrix}
1 & \frac{1}{\sqrt{2}} z_1 & \cdots & \frac{1}{\sqrt{2}} z_n
& -\sqrt{-1}
& \frac{\sqrt{-1}}{\sqrt{2}} z_1 & \cdots &
\frac{\sqrt{-1}}{\sqrt{2}} z_n\\
0 & \frac{1}{\sqrt{2}} & & & 0
&-\frac{\sqrt{-1}}{\sqrt{2}}\\
\vdots &&\ddots&&\vdots&&\ddots\\
0&&\cdots&\frac{1}{\sqrt{2}}& 0
&&\cdots&-\frac{\sqrt{-1}}{\sqrt{2}}
\end{pmatrix}.
\]
Since the Siegel manifold is complete and since the above
set is closed in $H_{n+1}$, the horizontal slice is
complete with respect to the Weil-Petersson metric.

\qed


\begin{thebibliography}{10}

\bibitem{AC}
D.~V. Alekseevsky and V.~Cort{\'e}s.
\newblock Classification of stationary compact homogeneous special
  pseudo-{K}\"ahler manifolds of semisimple groups.
\newblock {\em Proc. London Math. Soc. (3)}, 81(1):211--230, 2000.

\bibitem{BG}
R.~Bryant and P.~Griffiths.
\newblock Some observations on the infinitesimal period relations for regular
  threefolds with trivial canonical bundle.
\newblock In M.~Artin and J.~Tate, editors, {\em Arithmetic and Geometry},
  pages 77--85. Boston, Birkha\"user, 1983.


\bibitem{Gr}
P.~Griffiths, editor.
\newblock {\em Topics in Transcendental Algebraic Geometry}, volume 106 of {\em
  Ann. Math Studies}.
\newblock Princeton University Press, 1984.

\bibitem{Lu8}
Z.~Lu.
\newblock A note on special {K}\"ahler manifolds.
\newblock {\em Math. Ann.}, 313:711--713, 1999.

\bibitem{Lu5}
Z.~Lu.
\newblock On the {H}odge metric of the universal deformation space of
  {Calabi-Yau} threefolds.
\newblock {\em J. Geom. Anal.}, 11(1):103--118, 2001.

\bibitem{Lu20}
Z.~Lu.
\newblock On the rigidity of horizontal slices.
\newblock accepted by Portugaliae Mathematica, 2001.

\bibitem{S}
A.~Strominger.
\newblock {Special Geometry}.
\newblock {\em Comm. Math. Phy.}, 133:163--180, 1990.

\bibitem{T1}
G.~Tian.
\newblock Smoothness of the universal deformation space of compact {Calabi-Yau}
  manifolds and its {Peterson-Weil} metric.
\newblock In S.-T. Yau, editor, {\em Mathematical aspects of string theory},
  volume~1, pages 629--646. World Scientific, 1987.

\bibitem{Y1}
S.-T. Yau.
\newblock On the {R}icci curvature of a compact {K}\"ahler manifold and the
  complex {Monge-Ampere} equation, i,.
\newblock {\em Comm. Pure. Appl. Math}, 31:339--411, 1978.

\end{thebibliography}

\end{document}